\documentclass[11pt]{article}
\usepackage[numbers,sort&compress]{natbib}
\usepackage{graphicx,subfigure}
\usepackage{amsfonts}
\usepackage{mathrsfs}
\usepackage{amsmath}
\usepackage{amssymb}
\usepackage{graphicx}
\usepackage{xcolor}
\usepackage{epic}
\usepackage{amsthm,latexsym}
\usepackage{indentfirst}
\renewcommand{\paragraph}{\roman{paragraph}}
\usepackage[]{caption2} 
\newtheorem{theorem}{Theorem}[section]

\newtheorem{lemma}[theorem]{Lemma}

\begin {document}
\title{\bf  Maximizing Kirchhoff index of unicyclic graphs with fixed maximum degree}

\author{Dong Li$^{a}$, Xiang-Feng Pan$^{b,}$\thanks{Corresponding author. \textit{E-mail addresses:} djldong@126.com (D. Li), xfpan@ustc.edu (X.-F. Pan), liujiabaoad@163.com (J.-B. Liu), hqliu@hubu.edu.cn (H.-Q. Liu).},  Jia-Bao Liu$^{b}$, Hui-Qing Liu$^{a}$\\
{\small \it a.School of Mathematics and Statistics, Hubei University, Wuhan 430062, China}\\
{\small \it b.School of Mathematical Sciences, Anhui University, Hefei 230601, China}}
\date{}
\maketitle

\noindent {\bf Abstract.}  The Kirchhoff index of a connected graph is the sum of resistance distances
between all unordered pairs of vertices in the graph. Its considerable applications are found in
a variety of fields. In this paper, we determine the maximum value of Kirchhoff index among the
unicyclic graphs with fixed number of vertices and maximum degree, and characterize
the corresponding extremal graph.


\section {Introduction}
First we introduce some graph notations used in this paper. We only
consider finite, undirected and simple graphs. Other undefined
terminologies and notations may refer to \cite{Bon76}. Let $G$ be a connected graph with vertex set $V(G)$, whose vertices are labeled as $v_{1},v_{2},\ldots,v_{n}$. The famous Wiener index is considered as one of the most applicable graph invariant, used as one
of the topological indices for predicting physicochemical properties of organic compounds. The distance between vertices $v_{i}$ and $v_{j}$, denoted by $d(v_{i}, v_{j})$, is the length of a shortest path between them. The
Wiener index $W(G)$ \cite{west01} is the sum of distances between all pairs of vertices, that is,
$$ W(G)=\sum_{i<j}d(v_{i},v_{j}). $$

In 1993, Klein and Randic \cite{klein93} introduced a new distance function named resistance
distance basing on electrical network theory. They view graph $G$ as an electrical network such
that each edge of $G$ is assumed to be a unit resistor.  The resistance distance between
vertices $v_{i}$ and $v_{j}$, denoted by $r(v_{i}, v_{j})$ (if more than one graphs are considered, we use
$r_{G}(v_{i}, v_{j})$ to avoid confusion), is defined to be the effective resistance between nodes $v_{i}$ and
$v_{j}$ in $G$. Analogue to Wiener index, the Kirchhoff index $Kf(G)$ is defined as:
$$Kf(G)=\sum_{i<j}r(v_{i},v_{j}).$$
Let $Kf_{v_{i}}(G)$ denote the sum of all distances between vertex $v_{i}$ and the other vertices of
$G$, that is
$$Kf_{v_{i}}(G) = \sum _{j \neq i}r(v_{i}, v_{j}).$$

Kirchhoff index attracted extensive attention due to its wide applications in physics, chemistry, graph theory, etc. \cite{chen07,diudea01,Dobrynin01,entringer76,feng-yu14,fischermann03,fowler02,guo09,gutman86,liu15,liu14,liu-pan15,palacios01,qi13,trinajsti92,wiener1947,xiao04,zhang07,zhang09}. And some related results have been given in \cite{feng13,feng14,jelen03,li13,LiuPL2015,palacios01-1}.

 A graph $G$ is called a unicyclic graph
if it contains exactly one cycle. We may use the following notation \cite{yang08} to represent a unicyclic
graph:
$$G = U(C_{l}; T_{1}, T_{2},\ldots, T_{l}),$$
where $C_{l}$ is the unique cycle in $G$ with $V (C_{l}) = {v_{1}, v_{2}, \ldots , v_{l}}$ such that $v_{i}$ is adjacent to
$v_{i+1}$ (subscript module $l$) for $1 \leq i \leq l$. For each $i$, let $T_{i}$ be the component of $G - (V(C_{l}) - v_{i})$
containing $v_{i}$. Obviously $T_{i}$ is a tree. We say $T_{i}$ trivial if it is an isolated vertex. For example, see Fig. 1.

\begin{figure}[ht]
\centering
 \includegraphics[width=\textwidth]{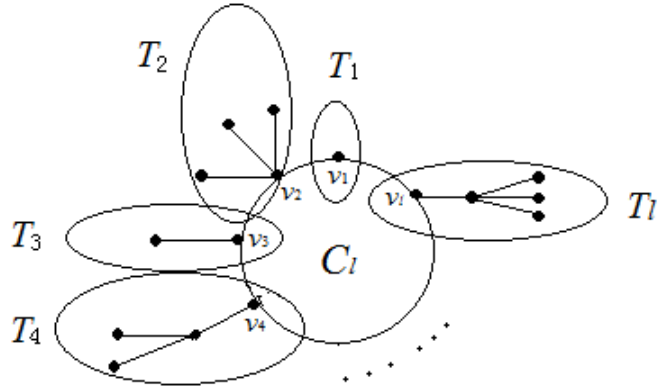}
\vspace{-5em} \caption{$U(C_{l}; T_{1}, T_{2}, \ldots, T_{l}).$}
\end{figure}

For convenience, we employ the following notation. Let $T_{n,\Delta}$ \cite{stevanovi08} be the tree on $n$ vertices obtained from a path $P_{n-\Delta+1}$ and an
empty graph $\overline{K}_{\Delta-1}$ by joining one end-vertex of $P_{n-\Delta+1}$ with every vertex of $\overline{K}_{\Delta-1}$.
Let $S_{n}^{l}$ \cite{yang08} denote the graph obtained from cycle $C_{l}$
by adding $n - l$ pendent edges to a vertex of $C_{l}$. Let $P_{n}^{l}$ \cite{yang08} denote the graph obtained by
identifying one end vertex of $P_{n - l + 1}$ with any vertex of $C_{l}$. And $\mathcal{P}_{n,\Delta}^{l}$ denote the set of the graphs obtained from $P_{t}^{l}(t\leq n-1)$ by attaching $n-t$ pendent edges to one vertex such that its degree achieves maximum degree $\Delta$ among all vertices. It is obvious that $S_{n}^{n}
 = P_{n}^{n} = C_{n}$.

In order to better
understand these graphs, we consider the set $\mathcal{G}_{n,l,\Delta}$, $\Delta \geq 2$, of connected graphs with
$n$ vertices containing cycle $C_{l}$ having the fixed value of the maximum degree $\Delta$. Otherwise, if we only bound
the maximum degree by $\Delta$, the graphs of maximum Kirchhoff index will inevitably be the paths. Still, even
with the requirement that a graph contains a vertex of degree $\Delta$, the extremal graphs
resemble a pathlike structure.

In this paper, we determine the maximum value of Kirchhoff index among the
unicyclic graphs with fixed number of vertices and maximum degree, and characterize
the extremal graphs. At the end of the article, we also give the conjecture on minimum Kirchhoff index among the
unicyclic graphs with fixed number of vertices and maximum degree, and
the extremal graphs.

\section {Main Result}
\begin{theorem}
For every graph $G \in \mathcal{G}_{n,l,\Delta}$, where $l\geq3$ and $n\geq l+\Delta-2$,
$$Kf(G)\le\frac{(n-\Delta+1)(n-\Delta+2)(n+2\Delta-9)}{6}+(\Delta-3)(\Delta-\frac{2}{3})+(n-\Delta)^{2}+\frac{7(n-\Delta)}{3}+2,$$ with equality holding
 if and only if $G\cong P_{n,\Delta}^{3}$ (see Fig. 2).
\end{theorem}

\begin{figure}[ht]
\centering
  \includegraphics[width=\textwidth]{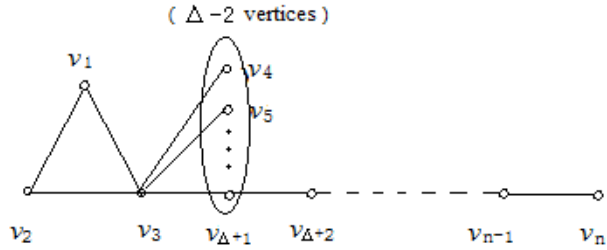}
\vspace{-10em} \caption{$P_{n,\Delta}^{3}$.}
\end{figure}

\section {The Proofs}
As we work with connected graphs only, the case $\Delta=1$
becomes possible only when $n=2$. Similarly, for $\Delta=2$ the only elements of $\mathcal{G}_{n,\Delta}$ are the
path $P_{n}$ and the cycle $C_{n}$. Thus, we may assume that $\Delta\geq3$ holds in the sequel. We start our proofs with a very
simple lemma.

\begin{lemma}\label{numedges} \cite{klein93}
 Let $x$ be a cutvertex of a graph, and let $a$ and $b$ be vertices occurring in
different components which arise upon deletion of $x$. Then
$$r(a, b) = r(a, x) + r(x, b).$$
\end{lemma}

\begin{lemma}\label{numedges} \cite{yang08}
Let $G = U(C_{l}; T_{1}, T_{2}, \ldots, T_{l})$. Then
$$Kf(G)=\sum_{i=1}^{l}W(T_{i})+\sum_{i<j}(|V(T_{j})|W_{v_{i}}(T_{i})+|V(T_{i})||V(T_{j})|\frac{(j-i)(l-j+i)}{l}+|V(T_{i})|W_{v_{j}}(T_{j})).$$
\end{lemma}

\begin{lemma}\label{numedges}  Assume that the vertex of maximum degree belongs to $V(T_{l})$. Let $G = U(C_{l}; T_{1}, T_{2}, \ldots, T_{l})$ and $G^{'} = U(C_{l}; P_{1}, P_{2}, \ldots, P_{l-1},T_{l})$. Then
$$Kf(G)<Kf(G^{'}).$$
\end{lemma}
{\bf Proof.} Since $W(T_{n})<W(P_{n})$, by the Lemma 3.2,
\begin{eqnarray}
  Kf(G)&=&\sum_{i=1}^{l}W(T_{i})+\sum_{i<j\leq l}(|V(T_{j})|W_{v_{i}}(T_{i})+|V(T_{i})||V(T_{j})|\frac{(j-i)(l-j+i)}{l}+|V(T_{i})|W_{v_{j}}(T_{j}))\nonumber
\\&<&\sum_{i=1}^{l}W(T_{i})+\sum_{i<j<l}(|V(P_{j})|W_{v_{i}}(P_{i})+|V(P_{i})||V(P_{j})|\frac{(j-i)(l-j+i)}{l}+|V(P_{i})|W_{v_{j}}(P_{j}))\nonumber
\\&&+|V(T_{l})|W_{v_{i}}(P_{i})+|V(P_{i})||V(T_{l})|\frac{(l-i)(0+i)}{l}+|V(P_{i})|W_{v_{l}}(T_{l})\nonumber
\\&=&Kf(G^{'}).\nonumber ~~~\square
\end{eqnarray}

\begin{lemma}\label{numedges}  Let $G\in\mathcal{G}_{n,l,\Delta}$, $G\notin \mathcal{P}_{n,\Delta}^{l}$ and assume that the vertex of maximum degree belongs to $V(T_{l})$. Then $Kf(G)<Kf(H)$ for every $H\in \mathcal{P}_{n,\Delta}^{l}.$
\end{lemma}
{\bf Proof.} Suppose that $G_{0} = U(C_{l}; T_{1}, T_{2}, \ldots, T_{l})$ has maximal Kirchhoff index among $\mathcal{G}_{n,l,\Delta}$.

{\bf Claim 1.} For each $i$, $T_{i}$ is a path with $v_{i}$ as one of its end vertices.

For each $i$, $W(T_{i})$ is maximal if and only if $T_{i}$ is a path and it is easy to
observe that $W_{v_{i}}(T_{i})$ is maximal if and only if $T_{i}$ is a path with $v_{i}$ as one of its end-vertices.
Hence Claim 1 holds by Lemma 3.2.

{\bf Claim 2.} If $l < n$, all but one of the $T_{i}$ are trivial.

Suppose to the contrary that there exists at least two trees such that they all have more
than one vertices. Let $T_{i}$ and $T_{j}$ be two such trees. By Claim 1, $T_{i}$ and $T_{j}$ are both paths.
Let $a\neq v_{i}$ and $b\neq v_{j}$ be endvertices of $T_{i}$ and $T_{j}$, respectively. Without losing generality,
assume that $Kf_{a}(G_{0}) \geq Kf_{b}(G_{0})$. Let $c$ be the neighbor of $b$ and $G_{0}^{'}=G_{0}-cb+ab$. Now
we show that $Kf(G_{0}^{'})> Kf(G_{0})$.

For any two vertices $x$, $y$ different from $b$, $r_{G_{0}}(x, y) = r_{G_{0}^{'}}(x, y)$, hence
$$Kf(G_{0})-Kf_{b}(G_{0}) = Kf(G_{0}^{'})-Kf_{b}(G_{0}^{'}).$$

On the other hand,
$$Kf_{b}(G_{0}^{'}) = Kf_{a}(G_{0}^{'}) + n-2 = Kf_{a}(G_{0}) + 1-r_{G_{0}}(a, b) + n-2.$$

Since $r_{G_{0}}(a, b) < d_{G_{0}}(a, b) < n-1$, it follows that $Kf_{b}(G_{0}^{'}) > Kf_{a}(G_{0}) \geq Kf_{b}(G_{0})$. Hence
$$Kf(G_{0}^{'}) = Kf(G_{0})-Kf_{b}(G_{0}) + Kf_{b}(G_{0}^{'}) > Kf(G_{0}).$$

This contradicts the choice of $G_{0}$, which implies Claim 2.

Claims 1 and 2 yield Lemma 3.4.~~~$\square$

\begin{lemma}\label{numedges} \cite{stevanovi08} For every graph $G\in\mathcal{G}_{n,\Delta}$, it holds that $$W(G)\leq W(T_{n,\Delta})$$
with equality if and only if $G$ is isomorphic to $T_{n,\Delta}$.
\end{lemma}

By the Lemma 3.3 and the Lemma 3.5, the number of the maximum degree less value is greater, i.e. when only one of the maximum degree value maximum. Therefore, to calculate the maximum, only need to discuss only one of the maximum degree situation can be.

\begin{lemma}\label{numedges}  Let $l\geq3$, then $\max\limits_{G\in \mathcal{P}_{n,\Delta}^{l}}\Big\{Kf(G)\Big\}$ achieves when the  maximum degree vertex is the connection point of the tree and the cycle.
\end{lemma}
{\bf Proof.} By Lemma 3.5, $Kf(G)\;(G\in\mathcal{P}_{n,\Delta}^{l})$ achieves maximum only if $G$ is one of two special graphs in $\mathcal{P}_{n,\Delta}^{l}$ (
see Fig. 3). Let the kirchhoff index of (a) and (b) in Fig. 3 be $Kf_{(a)}$ and $Kf_{(b)}$, respectively.

\begin{figure}[ht]
\centering
  \includegraphics[width=\textwidth]{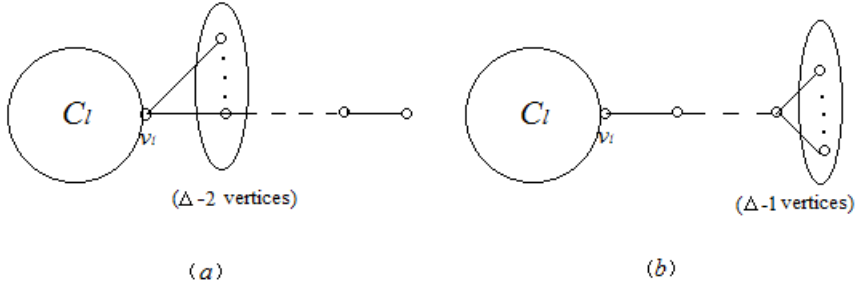}
\vspace{-10em} \caption{two special graphs in $\mathcal{P}_{n,\Delta}^{l}.$}
\end{figure}

By a routine calculation,
$$Kf(C_{l})=\frac{l^{3}-l}{12},$$

$$W(T_{n,\Delta}) =\Big(^{n-\Delta+2}_{~~~~3}\Big)+ (\Delta-1) \frac{(n-\Delta+1)(n-\Delta+2)}{2}+(\Delta-1)(\Delta-2)\nonumber $$
and
$$Kf_{v_{i}}(C_{l})=\frac{1\cdot(l-1)}{l}+\frac{2\cdot(l-2)}{l}+\cdots +\frac{(l-1)\cdot1}{l}=\frac{l^{2}-1}{6}.\nonumber $$
Then
\begin{eqnarray}
Kf_{(a)}&=&Kf(C_{l})+Kf(T_{n-l+1,\Delta-2})+(\Delta-2)Kf_{v_{l}}(C_{l})+(\Delta-2)(l-1)\cdot1\nonumber
\\&&+(n-l-\Delta+2)Kf_{v_{l}}(C_{l})+(l-1)\sum_{i=2}^{n-l-\Delta+3} i\nonumber
\\&=&\frac{l^{3}-l}{12}+\Big(^{n-l-\Delta+5}_{~~~~~3}\Big)+(\Delta-3)\frac{(n-l-\Delta+4)(n-l-\Delta+5)}{2}+(\Delta-3)(\Delta-4)\nonumber
\\&&+(\Delta-2)\frac{l^{2}-1}{6}+(\Delta-2)(l-1)+(n-l-\Delta+2)\frac{l^{2}-1}{6}\nonumber
\\&&+(l-1)\frac{(n-l-\Delta+2)(n-l-\Delta+5)}{2}\nonumber
\\&=&\frac{l^{3}-l}{12}+ \frac{(n-l-\Delta+5)(n-l-\Delta+4)(n-l+2\Delta-6)}{6}+\frac{(n-l)(l^{2}-1)}{6}\nonumber
\\&&+(\Delta-2)(l-1)+(\Delta-3)(\Delta-4)+\frac{(l-1)(n-l-\Delta+2)(n-l-\Delta+5)}{2} \nonumber
\end{eqnarray}
and
\begin{eqnarray}
Kf_{(b)}&=&Kf(C_{l})+Kf(T_{n-l+1,\Delta})+Kf_{v_{l}}(C_{l})(n-l-\Delta+1)+(l-1)\sum_{i=1}^{n-l-\Delta+1}i\nonumber
\\&&+(\Delta-1)Kf_{v_{l}}(C_{l})+(n-l-\Delta+2)(\Delta-1)(l-1)\nonumber
\\&=&\frac{l^{3}-l}{12}+\Big(^{n-\Delta-l+3}_{~~~~~3}\Big)+(\Delta-1)\frac{(n-l-\Delta+2)(n-l-\Delta+3)}{2}+(\Delta-1)(\Delta-2)+\nonumber
\\&&(n-l)\frac{l^{2}-1}{6}+(l-1)\frac{(n-l-\Delta+1)(n-l-\Delta+2)}{2}+(n-l-\Delta+2)(\Delta-1)(l-1)\nonumber
\\&=&\frac{l^{3}-l}{12}+\frac{(n-l-\Delta+2)(n-l-\Delta+3)(n-l+2\Delta-2)}{6}+\frac{(n-1)(l^{2}-1)}{6}+(\Delta-1)\nonumber
\\&&(\Delta-2)+(n-l-\Delta+2)(\Delta-1)(l-1)+\frac{(l-1)(n-l-\Delta+1)(n-l-\Delta+2)}{2}.\nonumber
\end{eqnarray}

Consequently,
\begin{eqnarray}
Kf_{(a)}-Kf_{(b)}&=&Kf(T_{n-l+1,\Delta-2})+(\Delta-2)(l-1)+(l-1)\frac{(n-l-\Delta+2)(n-l-\Delta+5)}{2}\nonumber
\\&&-\Big[Kf(T_{n-l+1,\Delta})+(n-l-\Delta+2)(\Delta-1)(l-1)\nonumber
\\&&+(l-1)\frac{(n-l-\Delta+1)(n-l-\Delta+2)}{2}\Big]\nonumber
\\&=&-5n+12\Delta+2n\Delta-2\Delta^{2}-18-2\Delta l+5l-4\Delta+10+(l-1)(-n\Delta+l\Delta\nonumber
\\&&+\Delta^{2}+3n-3l-4\Delta+4)\nonumber
\\&=&(\Delta-3)l^{2}+(\Delta^{2}-7\Delta-n\Delta+3n+12)l-8n+12\Delta+3n\Delta-3\Delta^{2}-12\nonumber
\\&=&(\Delta-3)l^{2}+(\Delta-3)(\Delta-n-4)l-8n+12\Delta+3n\Delta-3\Delta^{2}-12.\nonumber
\end{eqnarray}

When $\Delta=3$, we have $Kf_{(a)}-Kf_{(b)}=n-3>0$.
\\And when $\Delta>3$, we have
\begin{eqnarray}
Kf_{(a)}-Kf_{(b)}&=&(\Delta-3)l^{2}+(\Delta-3)(\Delta-n-4)l+3(\Delta-3)(n-\Delta+1)+(n-3)\nonumber
\\&=&l(\Delta-3)(l+\Delta-n-4)+3(\Delta-3)(n-\Delta+1)+(n-3)\nonumber
\\&>&3(\Delta-3)(l+\Delta-n-4)+3(\Delta-3)(n-\Delta+1)+(n-3)\nonumber
\\&=&3(\Delta-3)(l-3)+(n-3)\nonumber
\\&>&0.\nonumber
\end{eqnarray}

Hence Lemma 3.6 follows.~~~$\square$

\begin{lemma}\label{numedges} Let $G\in\mathcal{P}_{n,\Delta}^{l}, 3\leq l\leq n$, then $Kf(G)\leq Kf(P_{n,\Delta}^{3})$ and with equality holding
 if and only if $G\cong P_{n,\Delta}^{3}$.
\end{lemma}
{\bf Proof.} While the maximum degree vertex is the connection point of the tree and the cycle,
\begin{eqnarray}
Kf(P_{n,\Delta}^{3})&=&Kf(C_{3})+Kf(T_{n-2,\Delta-2})+\sum_{v_{i}\in C_{3}\setminus \{v_{3}\},v_{j}\in T_{n-2,\Delta-2}\setminus \{v_{3}\} } r_{v_{i}v_{j}}\nonumber
\\&=&2+\Big(^{n-\Delta+2}_{~~~~3}\Big)+(\Delta-3)\frac{(n-\Delta+1)(n-\Delta+2)}{2}+(\Delta-3)(\Delta-4)+\frac{10}{3}(\Delta-3)\nonumber
\\&&+2\Big[\frac{2}{3}(n-\Delta)+\frac{(n-\Delta)(n-\Delta+1)}{2}\Big]\nonumber
\\&=&\frac{(n-\Delta+1)(n-\Delta+2)(n+2\Delta-9)}{6}+(\Delta-3)(\Delta-\frac{2}{3})+(n-\Delta)^{2}+\frac{7}{3}(n-\Delta)\nonumber
\\&&+2.\nonumber
\end{eqnarray}

Hence we have
\begin{eqnarray}
Kf(P_{n,\Delta}^{3})-Kf_{(a)}&=&\frac{l^{3}}{12}-\frac{nl^{2}}{6}-(\Delta+\frac{1}{12}-2)l+n^{2}+\frac{5}{2}n-2n\Delta+\Delta^{2}+2\Delta-10-\frac{l^{3}}{3}+\nonumber
\\&&(\frac{n}{2}-\Delta+\frac{7}{2})l^{2}+(n\Delta-\frac{7}{2}n-\Delta^{2}+8\Delta-\frac{85}{6})l-n^{2}-n\Delta+5n+2\Delta^{2}\nonumber
\\&&-14n+22\nonumber
\\&=&-\frac{l^{3}}{4}+(\frac{7}{2}-\Delta)l^{2}+(-\Delta^{2}+7\Delta-\frac{49}{4})l+3\Delta^{2}-12\Delta+12+(\frac{l^{2}}{3}+l\Delta-\nonumber
\\&&\frac{7}{2}l-3\Delta+\frac{15}{2})n.\nonumber
\end{eqnarray}
And
$$\frac{l^{2}}{3}+l\Delta-\frac{7}{2}l-3\Delta+\frac{15}{2}>0,~~ n\geq l+\Delta-2.$$
So
\begin{eqnarray}
Kf(P_{n,\Delta}^{3})-Kf_{(a)}&\geq&-\frac{l^{3}}{4}+(\frac{7}{2}-\Delta)l^{2}+(-\Delta^{2}+7\Delta-\frac{49}{4})l+3\Delta^{2}-12\Delta+12\nonumber
\\&&+(\frac{l^{2}}{3}+l\Delta-\frac{7}{2}l-3\Delta+\frac{15}{2})(l+\Delta-2)\nonumber
\\&=&\frac{l^{3}}{12}+(\frac{1}{3}\Delta-\frac{2}{3})l^{2}+(-\frac{3}{2}\Delta+\frac{9}{4})l+\frac{3}{2}\Delta-3\nonumber
\\&\geq&0.\nonumber
\end{eqnarray}

Now we arrive at our main result.~~~$\square$

\section {Example}
When $n=100,~\Delta=96$, the maximum kirchhoff index is as follows, see Fig. 4.

\begin{figure}[ht]
\centering
  \includegraphics[width=\textwidth]{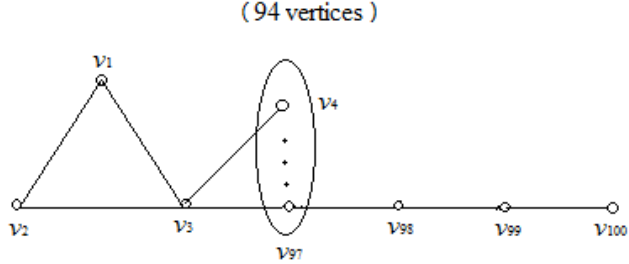}
\vspace{-12em} \caption{$P_{100,96}^{3}.$}
\end{figure}

\begin{eqnarray}
\max\limits_{G\in \mathcal{G}_{100,3,96}}\Big\{Kf(G)\Big\}&=&Kf(P_{100,96}^{3})\nonumber
\\&=&\frac{(100-96+1)(100-96+2)(100+2\times96-9)}{6}+(96-3)(96-\frac{2}{3})\nonumber
\\&&+(100-96)^{2}+\frac{7}{3}(100-96)+2\nonumber
\\&=&1415+8866+16+9\frac{1}{3}+2\nonumber
\\&=&10308\frac{1}{3}.\nonumber
\end{eqnarray}

\section {Conjecture}
Above the maximum Kirchhoff index among the
unicyclic graphs with fixed number of vertices and maximum degree is obtained. In the case of minimal values,  we end the text to the following conjecture:

(i) When $n\leq10$, or $n=11(\Delta\geq5)$, the Kirchhoff index among the
unicyclic graphs with fixed number of vertices and maximum degree has a minimum value as shown below, see Fig. 5.

\begin{figure}[ht]
\centering
  \includegraphics[width=\textwidth]{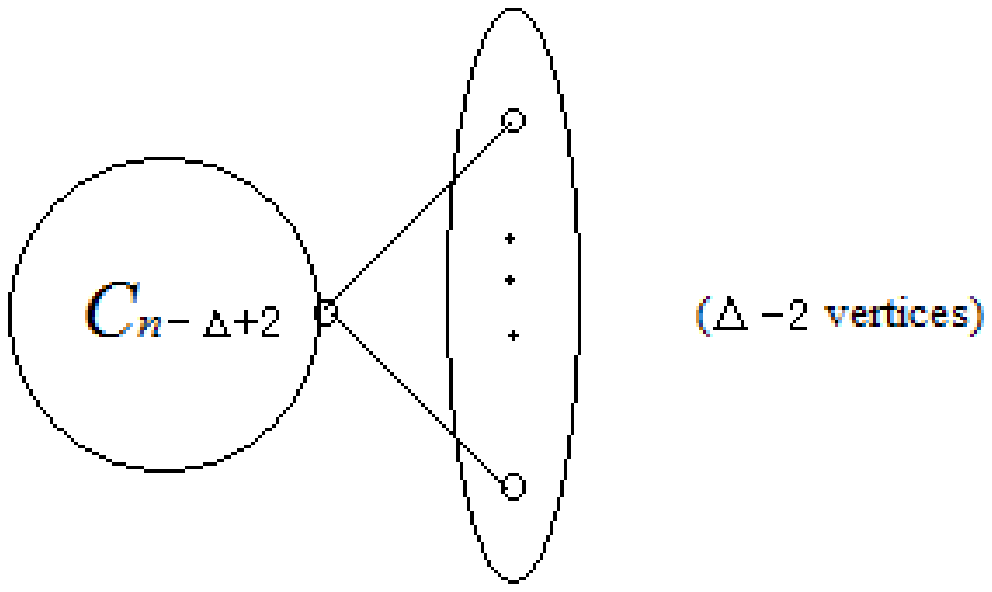}
\vspace{-8em} \caption{$G\in \mathcal{G}_{n,l,\Delta}$, with the minimum value of Kirchhoff index. }
\end{figure}

\begin{eqnarray}
\min\limits_{G\in \mathcal{G}_{n,l,\Delta}}\big\{Kf(G)\big\}&=&\frac{(n-\Delta+2)^{3}-(n-\Delta+2)}{12}+(\Delta-2)\Big[n-\Delta+2+\frac{(n-\Delta+2)^{2}-1}{6}\Big]\nonumber
\\&&+2\sum_{i=1}^{\Delta-3}i\nonumber
\\&=&\frac{(n+\Delta-2)(n-\Delta+3)(n-\Delta+1)}{12}+(\Delta-2)(n-1),\nonumber
\end{eqnarray}
where $l=n-\Delta+2$.

(ii) When $n\geq12$, or $n=11(\Delta\leq4)$, the Kirchhoff index among the
unicyclic graphs with fixed number of vertices and maximum degree has a minimum value as shown below, see Fig. 6.

\begin{figure}[ht]
\centering
  \includegraphics[width=\textwidth]{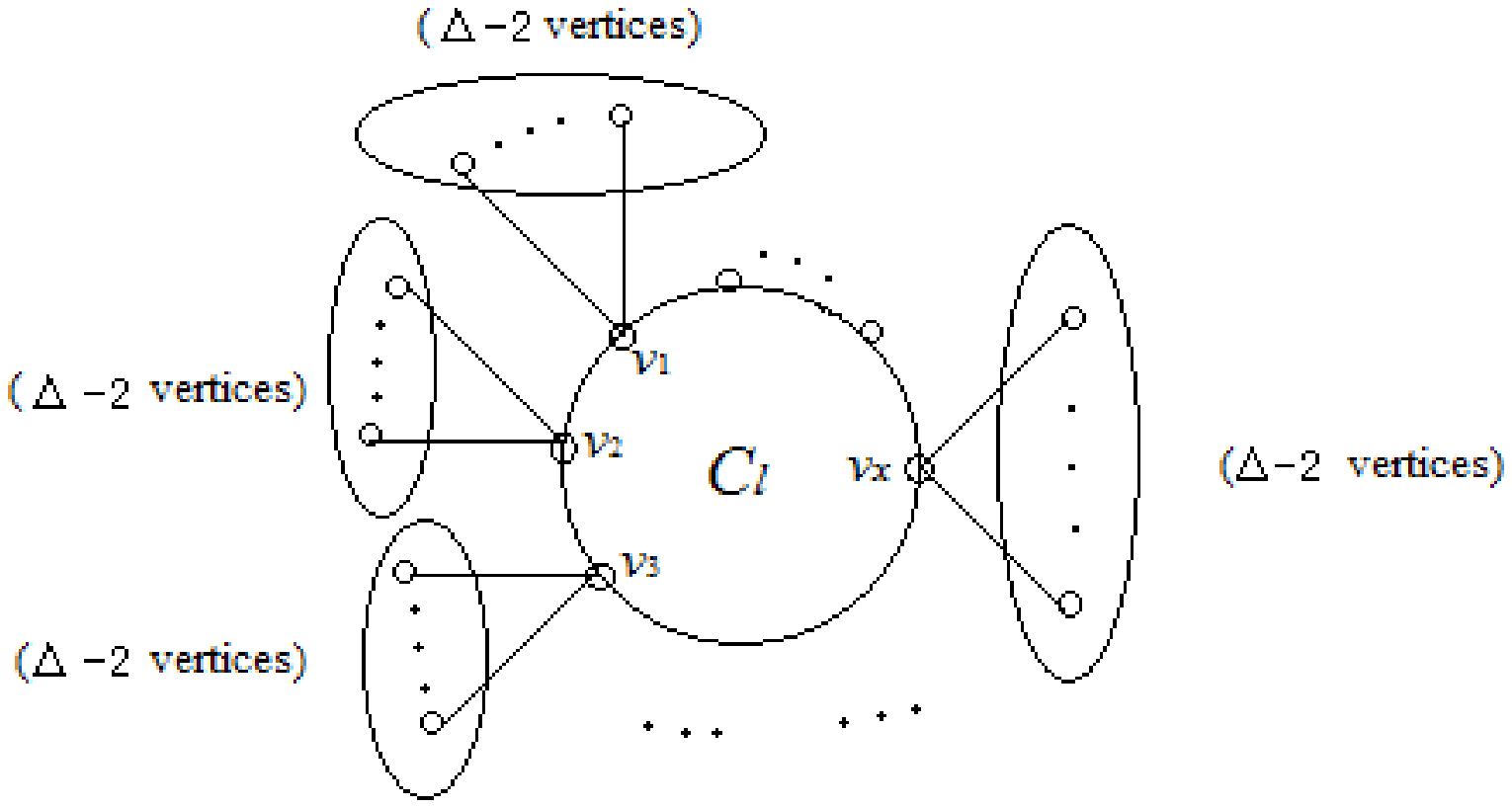}
\vspace{-10em} \caption{$G\in \mathcal{G}_{n,l,\Delta}$, with the minimum value of Kirchhoff index.  }
\end{figure}

\begin{eqnarray}
\min\limits_{G\in \mathcal{G}_{n,l,\Delta}}\big\{Kf(G)\big\}&=&Kf(C_{l})+x\cdot2\sum_{i=1}^{\Delta-3}i+x(\Delta-2)(\frac{l^{2}-1}{6}+l)\nonumber
\\&&+(x-1)(\Delta-2)^{2}\Big[2+\frac{(l-1)\cdot1}{l}\Big]+(x-2)(\Delta-2)^{2}\Big[2+\frac{(l-2)\cdot2}{l}\Big]\nonumber
\\&&+\cdots\nonumber
\\&&+1\cdot(\Delta-2)^{2}\Big[2+\frac{(l-x+1)\cdot(x-1)}{l}\Big]\nonumber
\\&=&\frac{l^{3}-l}{12}+x(\Delta-2)(\Delta-3)+x(\Delta-2)(\frac{l^{2}-1}{6}+l)+x(x-1)(\Delta-2)^{2}\nonumber
\\&&+(\Delta-2)^{2}\sum_{i=1}^{x-1}\frac{i(l-i)(x-i)}{l},\nonumber
\end{eqnarray}
where $l=n-x(\Delta-2)$, and $x\leq l\leq x+\Delta-2$, i.e. $\frac{n-\Delta+2}{\Delta-1}\leq x\leq \frac{n}{\Delta-1}$, $x$ is a positive integer.

\section*{Acknowledgments}
This work is partially supported by National Natural Science Foundation of
China (Nos. 11471016, 11401004, 11171097, 11371028 and 11571096), Natural
Science Foundation of Anhui Province of China (No. KJ2013B105),
Anhui Provincial Natural Science Foundation (No. 1408085QA03) and the Natural Science Foundation for the Higher Education Institutions of Anhui Province of China (No. KJ2015A331).

\end {document}